\documentclass[12pt,reqno]{amsart}
\usepackage[T1]{fontenc}
\usepackage[utf8]{inputenc}
\usepackage[a4paper,margin=3cm]{geometry}
\usepackage{latexsym,lmodern,color}
\usepackage[pagebackref, colorlinks=true, linkcolor=blue, citecolor=blue]{hyperref}
\usepackage{natbib}
\usepackage{comment}
\usepackage{amsmath,amsfonts,amssymb,amsthm}
\usepackage{stmaryrd}
\usepackage{graphicx}
\usepackage{esint}

\DeclareMathOperator{\N}{\mathbb{N}}
\DeclareMathOperator{\R}{\mathbb{R}}
\DeclareMathOperator{\C}{\mathbb{C}}
\DeclareMathOperator{\AAA}{\mathcal{A}}

\DeclareMathOperator{\CC}{\mathcal{C}}

\DeclareMathOperator{\1}{\mathbf{1}}
\DeclareMathOperator{\Reel}{\mathrm{Re}}
\DeclareMathOperator{\supp}{\mathrm{supp}}
\newcommand{\dd}{\mathrm{d}}

\newtheorem{Proposition}{Proposition}[section]
\newtheorem{Theorem}[Proposition]{Theorem}
\newtheorem*{Theorem*}{Theorem}
\newtheorem{Lemma}[Proposition]{Lemma}

\theoremstyle{definition}

\newtheorem*{Remark}{Remark}

\newtheorem*{Proof}{Proof}

\title[Free Fokker-Planck equation]{Convergence to Equilibrium in the Free Fokker-Planck Equation With a Double-Well Potential}
\author{Catherine \textsc{Donati-Martin}} 
\address[CDM]{Laboratoire de Math\'ematiques de Versailles, UVSQ, CNRS, Universit\'e Paris-Saclay, 45 avenue des \'Etats-Unis, 78035 \textsc{Versailles Cedex}, France. E-mail: \texttt{catherine.donati-martin@uvsq.fr} } 
\author{Benjamin \textsc{Groux}}
\address[BG]{Laboratoire de Math\'ematiques de Versailles, \ UVSQ, \ CNRS, \ Universit\'e Paris-Saclay, 45 avenue des \'Etats-Unis, 78035 \textsc{Versailles Cedex}, France. E-mail: \texttt{benjamin.groux@uvsq.fr} } 
\author{Myl\`ene \textsc{Ma\"ida}}
\address[MM]{Universit\'e Lille 1, Laboratoire Paul Painlev\'e, Cit\'e Scientifique, 59655 \textsc{Villeneuve d'Ascq Cedex}, France. E-mail: \texttt{mylene.maida@math.univ-lille1.fr} }   
\date{\today}

\begin{document}

\maketitle

\begin{abstract}
We consider the one-dimensional free Fokker-Planck equation
\begin{displaymath}
\frac{\partial \mu_t}{\partial t} = \frac{\partial}{\partial x} \left[ \mu_t \cdot \left( \frac12 V' - H\mu_t \right) \right] \, ,
\end{displaymath}
where $H$ denotes the Hilbert transform and $V$ is a particular double-well quartic potential, namely $V(x) = \frac14 x^4 + \frac{c}{2} x^2$, with $c \ge -2$. We prove that the solution $(\mu_t)_{t \ge 0}$ of this PDE converges in Wasserstein distance of any order $p \ge 1$ to the equilibrium measure $\mu_V$ as $t$ goes to infinity. This provides a first result of convergence for this equation in a non-convex setting. The proof involves free probability and complex analysis techniques.\\
\end{abstract}

\textbf{AMS 2010 Classification Subject.} 35B40, 46L54, 60B20.\\


\section{Introduction}

We consider the following \emph{one-dimensional free Fokker-Planck equation}
\begin{equation} \label{FP}
\frac{\partial \mu_t}{\partial t} = \frac{\partial}{\partial x} \left[ \mu_t \cdot \left( \frac12 V' - H\mu_t \right) \right] \, .
\end{equation}
In this equation, $(\mu_t)_{t \ge 0}$ denotes a family of probability measures on $\R$, $V : \R \to \R$ is a given potential, and $H$ denotes the Hilbert transform, that is, for any probability measure $\mu$ on $\R$ and $x \in \R$,
\begin{displaymath}
H\mu(x) = \fint_{\R} \frac{1}{x-y} \, \dd\mu(y) := \lim_{\epsilon \downarrow 0} \int_{\R \setminus[x-\epsilon, x+\epsilon]} \frac{1}{x-y} \, \dd\mu(y) \, ,
\end{displaymath}
where $\fint$ stands for the principal value of the integral. Partial differential equation (PDE) \eqref{FP} must be understood in the sense of distributions, i.e. for any regular enough test function $\varphi : \R \to \R$,
\begin{displaymath}
\frac{\dd}{\dd t} \int \varphi(x) \, \dd\mu_t(x) = - \frac12 \int V'(x)\varphi'(x) \, \dd\mu_t(x) + \frac12 \iint \frac{\varphi'(x)-\varphi'(y)}{x-y} \, \dd\mu_t(x) \dd\mu_t(y) \, .
\end{displaymath}
Under this form, it is sometimes called the \emph{McKean-Vlasov equation with logarithmic interaction}.

\subsection{Existence and uniqueness}

As far as we know, the problems of existence and uniqueness of the solution to the PDE \eqref{FP} are not completely solved.

Existence was tackled in \cite[Theorem 3.1]{BS}. Using the free stochastic calculus formalism (see \citep{Biane, BSp1} for an introduction), they proved that if, roughly speaking, $V$ is locally Lipschitz and grows "nicely" at infinity then, for any initial condition $X_0$ whose distribution is compactly supported, the free stochastic differential equation (SDE)
\begin{equation} \label{freeSDE}
\dd X_t = \dd S_t - \frac12 V'(X_t) \dd t \, ,
\end{equation}
where $S$ is a free Brownian motion, admits a unique solution $(X_t)_{t \ge 0}$ starting from $X_0$.
As they also checked that the distribution of the solution $(X_t)_{t \ge 0}$ satisfies the PDE \eqref{FP}, this proves the existence of a solution for the latter.

As for uniqueness, using free transportation techniques, \cite[Theorem 1.3]{LLX} shows that the free Fokker-Planck equation \eqref{FP} admits a unique solution starting from a compactly supported $\mu_0$ as soon as $V$ satisfies the same properties as in \cite[Theorem 3.1]{BS} and $V''$ is uniformly bounded below.\\

In this paper, we are interested in the free Fokker-Planck equation \eqref{FP} for the particular potential
\begin{equation} \label{V}
V(x) = \frac14 x^4 + \frac{c}{2} x^2 \, , \qquad c \ge -2 \, .
\end{equation}

\begin{figure}[ht]
\centering
\includegraphics[scale=1]{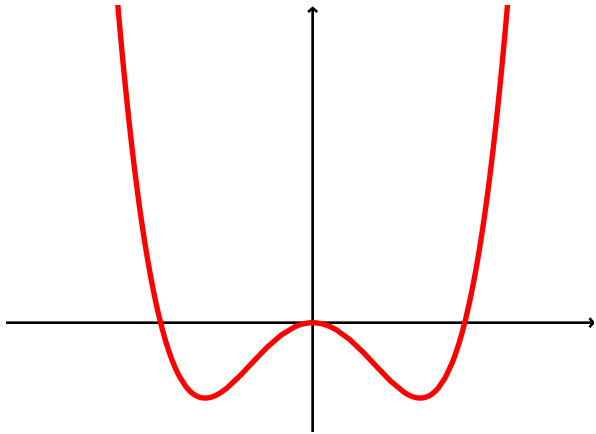}
\caption{Potential $V$ defined by \eqref{V} with $-2 \le c < 0$.}
\end{figure}
Indeed, the quadratic potential $V(x) = \frac{x^2}{2}$ gives rise to the \emph{free Ornstein-Uhlenbeck process} and is well understood; the quartic potential \eqref{V} is then the most simple example of a potential satisfying the assumptions of \cite[Theorem 3.1]{BS} and \cite[Theorem 1.3]{LLX}. Consequently, given any compactly supported initial condition, existence and uniqueness of the solution to Equation \eqref{FP} are ensured and we can moreover identify this solution with the distribution of the solution $(X_t)_{t \ge 0}$ to the free SDE \eqref{freeSDE}.\\

Note that the case when $c \ge 0$ is already covered in the existing literature, and the aim of this paper is to extend the study of the asymptotic behaviour of the solution $(\mu_t)_{t \ge 0}$ to the range  $c \ge -2$ (see also \cite[Conjecture 5.1]{LLX}).

\subsection{Granular media equation} \label{section_granular}
We say that a family $(\mu_t)_{t \ge 0}$ of probability measures on $\R^d$ having densities $(\rho_t)_{t \ge 0}$ satisfies a \emph{granular media equation} if we have in distribution
\begin{equation} \label{granular_eq}
\frac{\partial \mu_t}{\partial t} = \nabla \cdot \left[ \mu_t \nabla (\mathcal{U}'(\rho_t) + \mathcal{V} + \mathcal{W}*\rho_t) \right] \, ,
\end{equation}
where $\mathcal{U} : \R_+ \to \R$ can be seen as an internal energy, $\mathcal{V} : \R^d \to \R$ as a confinement potential, $\mathcal{W} : \R^d \to \R$ as an interaction potential, and where the operation $*$ is the usual convolution in $\R^d$.

The free Fokker-Planck equation \eqref{FP} corresponds to the particular case
\begin{displaymath}
d=1, \ \mathcal{U}(s)=0, \ \mathcal{V}(x)=\frac12 V(x), \ \textrm{and} \ \mathcal{W}(x)=-\log|x| \, .
\end{displaymath}

Several classical partial differential equations arising from physics fall into the class of a granular media equation (see \cite[Chapter 9.1]{Vill}), starting from the heat equation (for $\mathcal{U}(s)=s\log(s)$, $\mathcal{V}=0$, $\mathcal{W}=0$).  Conditions are known to ensure that Equation \eqref{granular_eq} admits a unique solution, but we will not discuss this point here and rather focus on reviewing some existing results about the long-time behaviour of the solutions. 

At least two main techniques can be identified. We focus on the \emph{entropy dissipation method}, which is close to the techniques we will use in this work, but we will also briefly mention some results using an \emph{approximation by a particle system}.\\

The first results using entropy dissipation are due to \cite{BCP, BCCP}, who are interested in particular potentials arising from physics. In these works, to study the long-time behaviour of the solution, the authors consider an entropy functional naturally associated to \eqref{granular_eq} defined by
\begin{displaymath}
F(\mu) = \int_{\R^d} \mathcal{U}(\mu(x)) \, \dd x + \int_{\R^d} \mathcal{V}(x) \, \dd\mu(x) + \frac12 \iint_{\R^d \times \R^d} \mathcal{W}(x-y) \, \dd\mu(x) \dd\mu(y) \, ,
\end{displaymath}
as the sum of an internal energy, a potential energy, and an interaction energy associated to a given measure $\mu$, and they show that this entropy is strictly decreasing along the trajectory $(\mu_t)_{t \ge 0}$. Under appropriate assumptions on $\mathcal{V}$ and $\mathcal{W}$, $F$ admits a unique minimizer $\mu_{\infty}$, which is shown to be the limit of $\mu_t$ as $t \to +\infty$.

Combining this entropy dissipation method with optimal transport techniques, \cite{CMCV} establish convergence of $\mu_t$ in a more general setting, even leading in some cases to explicit rates of convergence. Note that \cite{CGM, BGM, BGG1, BGG2} got various improvements of these results. Nevertheless, all these results require convexity, positivity, or smoothness assumptions on $\mathcal{V}$ and $\mathcal{W}$.

In a series of works, \cite{Tugaut1, Tugaut2} then tackled the problem of non-convex potentials in the case when $\mathcal{U}(s) = \sigma s\log(s)$ for a small $\sigma>0$ but his results still require a smooth interaction $\mathcal{W}$.\\

An example of physically meaningful singular interaction is given by $\mathcal{W}(x) = -\log|x|$, which is out of reach of the previous methods as they are. This problem of the logarithmic interaction in dimension one has been recently tackled by \cite{LLX}, who could adapt Carrillo, McCann, and Villani's method to the free probability framework, at least in the case of a convex potential $\mathcal{V}$.
We also mention the recent work by \cite{CCV}, in which a two-dimensional logarithmic interaction is considered, corresponding to the Keller-Segel model.

In view of these results, the study of the long-time behaviour of the solution of the granular media equation with a logarithmic interaction and a non-convex potential, such as \eqref{V}, is of natural interest.\\

Another motivation to study the granular media equation with a logarithmic interaction is its link with a particle system which is well known in random matrix theory.

For any $N \ge 1$, we consider the system of stochastic differential equations (SDEs)
\begin{equation} \label{system_SDEs}
\forall i \in \llbracket 1,N \rrbracket , \ \dd X_i^N(t) = \sqrt{2} \dd B_i(t) - \nabla \mathcal{V}(X_i^N(t)) \dd t - \frac1N \sum_{j \neq i} \nabla \mathcal{W}(X_i^N(t) - X_j^N(t)) \dd t \, ,
\end{equation}
where the $B_i$'s are independent Brownian motions. The solution $(X_1^N(t), \ldots X_N^N(t))_{t \ge 0}$ of \eqref{system_SDEs} is a natural particle system that can be associated to PDE \eqref{granular_eq}.

Indeed, as the number of particles $N$ goes to infinity, its empirical measure $\left(\frac1N \sum_{i=1}^N \delta_{X_i^N(t)}\right)_{t \ge 0}$ converges to a solution of PDE \eqref{granular_eq}, the Brownian term giving rise to an internal energy $\mathcal{U}(s) = s\log(s)$. If the Brownian term in \eqref{system_SDEs} is multiplied by $\frac{1}{\sqrt{N}}$, it disappears in the limit and we get the solution of \eqref{granular_eq} with $\mathcal{U}(s)=0$.

Using propagation of chaos for a modified approximating particle system, \cite{Malrieu} recovered some of the results of \citep{CMCV}. \cite{BGM} and \cite{CGM} also considered a particle system to prove convergence of the solution to the PDE they study.

In the case of a logarithmic interaction ($\mathcal{U}(s) = 0$, $\mathcal{V}(x) = 0$, and $\mathcal{W}(x) = -\log|x|$), the particle system \eqref{system_SDEs} is the well-known \emph{Dyson Brownian motion} introduced by \cite{Dyson} as the process of eigenvalues of Hermitian random matrices with Brownian entries. This process (or its variant when $\mathcal{V}$ is quadratic) has been much studied, among others by \cite{Chan, CL, Font, RS}; see also \cite[Section 4.3]{AGZ}. For more general potentials $\mathcal{V}$, this particle system has been studied by \cite{AD} in the cubic case and by \cite{LLX} in the convex case.

With this point of view, the reader who is familiar with random matrix theory will get an insight why a natural candidate for the long-time limit of the solution to the free Fokker-Planck equation \eqref{FP} should be the equilibrium measure associated to potential $V$, that we now define.

\subsection{Main result of the paper} \label{section_result}

Let $D$ be a closed subset of $\C$ and $V : D \to \C$ be a polynomial such that
\begin{displaymath}
\lim_{|z| \to +\infty, \, z \in D} \Reel V(z) - 2\log|z| = +\infty \, .
\end{displaymath}
Then the functional
\begin{equation} \label{SigmaV}
\Sigma_V : \mu \mapsto - \iint_{D^2} \log|z-t| \, \dd\mu(z) \dd\mu(t) + \int_{D} \Reel V(z) \, \dd\mu(z) \, ,
\end{equation}
called \emph{Voiculescu free entropy}, admits a unique minimizer among probability measures supported on $D$. This minimizer is called the \emph{equilibrium measure} associated to $V$ and $D$, and is denoted by $\mu_V$. Note that when $D \subset \R$ and $V$ is real-valued, we have
\begin{displaymath}
\Sigma_V(\mu) = - \iint_{\R^2} \log|x-y| \, \dd\mu(x) \dd\mu(y) + \int_{\R} V(x) \, \dd\mu(x) \, .
\end{displaymath}
We refer to \citep{ST} for a development on this topic for which the equilibrium measure is defined in a much more general setting.

For the quartic potential
\begin{displaymath}
V(x) = \frac14 x^4 + \frac{c}{2} x^2
\end{displaymath}
and $D = \R$, the equilibrium measure is explicitly known (see \cite[Example 3.2]{J} for instance):
\begin{itemize}
\item when $c \ge -2$, its density is given by
\begin{equation} \label{mu_V_c>-2}
\rho_V(x) = \frac{1}{\pi} \left( \frac12 x^2 + b_0 \right) \sqrt{a^2-x^2} \1_{[-a,a]}(x)
\end{equation}
where
\begin{displaymath}
a^2 = \frac23 \left( \sqrt{c^2+12} - c \right) , \qquad b_0 = \frac13 \left( c + \sqrt{\frac{c^2}{4} + 3} \right) \, ,
\end{displaymath}
\item when $c < -2$, its density is given by
\begin{equation} \label{mu_V_c<-2}
\rho_V(x) = \frac{1}{2\pi} |x| \sqrt{(x^2-a^2)(b^2-x^2)} \1_{[-b,-a] \cup [a,b]}(x)
\end{equation}
where $a^2 = -2-c$, $b^2 = 2-c$.\\
\end{itemize}

In this paper, we focus on the case when $c \ge -2$, in which the equilibrium measure has connected support. Here is the main result of this paper.

For any real $p \ge 1$, if $\mu$ and $\nu$ are probability measures on $\R$ such that $|\cdot |^p$ is integrable for $\mu$ and $\nu$, then the \emph{Wasserstein distance} of order $p$ between $\mu$ and $\nu$ is defined by
\begin{displaymath}
W_p(\mu,\nu) = \Bigr(\inf\iint |x-y|^p\pi(\dd x,\dd y)\Bigr)^{1/p} \, ,
\end{displaymath}
where the infimum runs over all probability measures $\pi$ on $\R\times\R$ with marginal distributions $\mu$ and $\nu$.

\begin{Theorem} \label{theorem}
Let $V(x) = \frac14 x^4 + \frac{c}{2} x^2$ with $c \ge -2$. Given any compactly supported probability measure $\mu_0$ on $\R$, the solution $(\mu_t)_{t \ge 0}$ of the free Fokker-Planck equation \eqref{FP} with initial condition $\mu_0$ satisfies
\begin{displaymath}
\lim_{t \to +\infty} W_p(\mu_t,\mu_V) = 0
\end{displaymath}
for all $p \ge 1$, where $\mu_V$ is given by \eqref{mu_V_c>-2}.
\end{Theorem}

Let us mention that the convergence in $W_p$ distance is equivalent to the weak convergence of measures together with the convergence of the moments of order $p$ (see for instance \cite[Theorem 7.12]{Vill}).\\

The case when $c \ge 0$ was already covered by previous results of \cite{LLX}. Indeed, in \cite[Theorem 1.6 (i) and (ii)]{LLX}, they proved that
\begin{displaymath}
\lim_{t \to +\infty} W_2(\mu_t,\mu_V) = 0
\end{displaymath}
as soon as $V$ is convex (in the case when $V$ is strictly convex, they even get that $t \mapsto W_2(\mu_t,\mu_V)$ exponentially decreases to 0). They provide a proof for convergence in $W_p$, $p \le 2$, that could be easily extended to any $W_p$, $p>2$. We also refer to \cite[Section 6.4]{G} for some complements.\\

On the other hand, a result of \cite[Section 7.1]{BS} implies that, if $c<0$ and $|c|$ is large enough, then there exist initial conditions $\mu_0$ for which the solution $(\mu_t)_{t \ge 0}$ does not converge towards the equilibrium measure $\mu_V$.\\

The rest of the paper is organized as follows. Some tools, such as properties of $(\mu_t)_{t \ge 0}$ viewed as the law of a free diffusion and the description of critical measures via complex analysis techniques, are introduced in Section \ref{section_tools}, and Section \ref{section_proof} uses these tools to prove Theorem \ref{theorem}. Section \ref{section_perspectives} is the final section of this paper, in which we present some perspectives for future work.

\section{Free probability and complex analysis tools} \label{section_tools}

\subsection{Some properties of the solution of the free Fokker-Planck equation} \label{section_properties_freeSDE}

As we explained in the introduction of the paper, the solution of the free Fokker-Planck equation \eqref{FP} can be interpreted as the distribution of the solution to the free SDE
\begin{displaymath}
\dd X_t = \dd S_t - \frac12 V'(X_t) \dd t \, ,
\end{displaymath}
where $S$ is a free Brownian motion. As a consequence, it inherits some properties of free diffusions with regular drift, studied by \cite{BS}, the most important of which are the following.

\begin{Proposition}[see {\cite[Theorems 3.1 and 5.2]{BS}}] \label{lemma1}
Let $V$ be a $\CC^1$ potential such that $V'$ is locally Lipschitz, and such that there exist $a<0$ and $b>0$ such that, for all $x \in \R$,
\begin{equation} \label{condition_BSp}
-xV'(x) \le ax^2+b \, .
\end{equation}
Let $(\mu_t)_{t \ge 0}$ be the solution of the free Fokker-Planck equation \eqref{FP} starting from a compactly supported $\mu_0$.
\begin{itemize}
\item[(i)] There exists $M>0$ such that, for every $t>0$,
\begin{equation} \label{BianeSpeicher1}
\supp(\mu_t) \subset [-M,M] \, .
\end{equation}
\item[(ii)] There exist $K_1,K_2>0$ depending only on $V$ such that, for every $t>0$, the density $\rho_t$ of $\mu_t$ satisfies
\begin{equation} \label{BianeSpeicher2}
\| \rho_t \|_{\infty} \le \frac{K_1}{\sqrt{t}} + K_2 , \qquad \| D^{1/2} \rho_t \|_2 \le \frac{K_1}{t} + K_2 \, ,
\end{equation}
where $D^{1/2}$ is the fractional derivative of order $1/2$.
\item[(iii)] The family $\{ \rho_t \}_{t \ge 1}$ lives in a subset $\AAA$ of $L^2([-M,M])$ which is compact for the topology induced by $\| \cdot \|_2$.
\end{itemize}
\end{Proposition}

In the statement of Point (ii), the notion of half-derivative appears. It can be defined by several ways; we will just mention that for $u \in L^2$, the derivative of order $1/2$ of $u$ is the inverse Fourier transform of $\xi \mapsto (1+\xi^2)^{1/4} \hat{u}(\xi)$, where $\hat{u}$ is the Fourier transform of $u$ (see \cite[Chapter 4]{Demengel} for instance).\\

We notice that the potential \eqref{V} satisfies the assumptions of Proposition \ref{lemma1}. Nevertheless, we include here a proof of Point (iii), since it will play a key role in extending the entropy dissipation method to this singular interaction, as we will explain at the beginning of Section \ref{section_proof}.

\begin{Proof}[of Proposition \ref{lemma1} (iii)]
By Proposition \ref{lemma1} (i)-(ii), there exist $M, K_1, K_2 > 0$ such that, for every $t>0$, \eqref{BianeSpeicher1} and \eqref{BianeSpeicher2} hold. For every $t>0$, we denote by $\AAA_t$ the set of probability density functions $f$ with support in $[-M,M]$ which satisfy $\|f\|_{\infty} \le \frac{K_1}{\sqrt{t}} + K_2$ and $\| D^{1/2}f \|_2 \le \frac{K_1}{t} + K_2$. Note that, for $t>0$, $\AAA_t$ contains all the $\rho_{t+s}$'s, $s \ge 0$, where $\rho_{t+s}$ denotes the density of the measure $\mu_{t+s}$ as in Point (ii).\\

Furthermore, for every $t>0$, $\AAA_t$ is a subset of the Sobolev space $H^{1/2}([-M,M])$, defined as the set of $L^2$-probability density functions whose derivative of order $1/2$ belongs to $L^2$. Because the injection of $H^{1/2}([-M,M])$ in $L^p([-M,M])$ is compact for every $p \in [1,\infty)$ (see \cite[Theorem 4.54]{Demengel} for instance) and $\AAA_t$ is bounded in $H^{1/2}([-M,M])$, we can deduce that the set $\AAA_t$ is relatively compact in $L^2([-M,M])$. Hence, we can choose for $\AAA$ the closure of $\AAA_1$ in $L^2([-M,M])$. \qed
\end{Proof}

Furthermore, even if we consider here a singular interaction in granular media equation \eqref{granular_eq}, thanks to the bounds stated in Proposition \ref{lemma1}, we have an entropy dissipation formula as in \citep{CMCV}.

\begin{Proposition}[see {\cite[Proposition 6.1]{BS}}]
Under the assumptions and notations of Proposition \ref{lemma1}, we have
\begin{equation} \label{derivativeSigmaV}
\frac{\dd}{\dd t} \Sigma_V(\mu_t) = -2 \int \left| \frac12 V' - H\mu_t \right|^2 \, \dd\mu_t \, .
\end{equation}
\end{Proposition}

As this formula suggests, probability measures $\mu$ supported in $\mathbb{R}$ satisfying the \emph{Euler-Lagrange equation}
\begin{equation} \label{EulerLagrange}
H\mu = \frac12 V' \quad \mu\mathrm{-a.e.}
\end{equation}
for a real potential $V$ are the candidates for the long-time limit of $\mu_t$. We will call these measures \emph{stationary measures} because they are exactly the stationary solutions of Equation \eqref{FP} in the sense of PDEs, i.e. they are solutions of Equation \eqref{FP} that are constant in time.

\subsection{Critical measures and their identification} \label{section_critical}

In addition to the equilibrium measure and stationary measures that we encounter in our problem, we introduce the notion of critical measure, as defined by \cite{MFR}.

A probability measure $\mu$ on $\C$ such that $\Sigma_V(\mu) < +\infty$ is called a \emph{critical measure} associated to $V$ if, for every $h : \C \to \C$ regular enough, the quantity
\begin{displaymath}
D_h\Sigma_V(\mu) = \lim_{s \to 0} \frac{\Sigma_V(\mu^{h,s}) - \Sigma_V(\mu)}{s}
\end{displaymath}
is zero, where $\mu^{h,s}$ is the push-forward measure of $\mu$ by the deformation of identity $z \mapsto z + sh(z)$, $s \in \C$.\\

The reason why we consider critical measures is that, by \cite[Lemma 3.7]{MFR}, we have
\begin{displaymath}
D_h\Sigma_V(\mu) = \Reel \left( \int V'(x)h(x) \, \dd\mu(x) - \iint \frac{h(x)-h(y)}{x-y} \, \dd\mu(x) \dd\mu(y) \right) \, ,
\end{displaymath}
hence for a probability measure $\mu$ supported on $\R$, the previous condition is equivalent to the Euler-Lagrange equation \eqref{EulerLagrange}. As a result, critical measures supported on $\R$ are exactly stationary measures, and we will be able to use some tools developed to identify critical measures in order to identify stationary measures.

Note that in general, several critical measures may exist while there is only one equilibrium measure. This is the case for a potential satisfying the conditions given in \cite[Section 7.1]{BS} for instance. A key point in the proof of Theorem \ref{theorem} will be to show that for the quartic potential \eqref{V}, there is no other critical measure than the equilibrium measure.\\

The following statement gives the most important properties of  critical measures supported in $\mathbb R$ we will use in the sequel. The key point is that the Stieltjes transform of a critical measure $\mu$, defined on $\C \setminus \R$ by
\begin{displaymath}
G_{\mu}(z) = \int_{\R} \frac{1}{z-x} \, \dd\mu(x) \, ,
\end{displaymath}
satisfies an algebraic equation.

\begin{Proposition}[see \cite{KS, HKL}] \label{lemma2}
Let $V$ be a polynomial and $\mu$ be a critical measure supported on $\R$.
\begin{itemize}
\item[(i)] There exists a polynomial $R$ of degree $2\deg(V)-2$ such that
\begin{equation} \label{R_def}
R(z) = \left( \frac12 V'(z) - G_{\mu}(z) \right)^2
\end{equation}
almost everywhere for Lebesgue measure on $\C$. Moreover, we have
\begin{equation} \label{R_expr}
R(z) = \frac14 V'(z)^2 - \int_{\R} \frac{V'(x)-V'(z)}{x-z} \, \dd\mu(x) \, .
\end{equation}
\item[(ii)] Every non-real root of $R$ has even multiplicity.
\item[(iii)] The support of $\mu$ is a finite union of intervals connecting zeros of $R$.
\end{itemize}
\end{Proposition}

Point (i) combines Proposition 3.7 and Formula (3.31) from \cite{KS}. Point (ii) is an easy consequence of analyticity of Stieltjes transform, see \cite[Lemma 2.6]{HKL}. At last, Point (iii) comes from \cite[Proposition 3.9]{KS}.

Let us remark that a critical measure $\mu$ is completely determined by the associated polynomial $R$. Hence, finding critical measures boils down to determining all possible polynomials satisfying Equations \eqref{R_def} and \eqref{R_expr}. For the quartic potential and other polynomials with few monics, this is possible to do so (see \cite{MFR, HKL} for examples). However, in the quartic case, we will only use $R$ in order to show that a critical measure has connected support. Indeed, as soon as this is the case, we can just recover $\mu$ by solving a singular integral equation, as we will do in the next subsection.

\begin{Proposition} \label{lemma2_quartic}
For the potential $V(x) = \frac14 x^4 + \frac{c}{2} x^2$ with $c \ge -2$, every critical measure supported in $\mathbb{R}$ has connected support.
\end{Proposition}

\begin{Proof}
By \eqref{R_expr}, the polynomial $R$ defined in \eqref{R_def} is given by
\begin{displaymath}
R(z) = \frac14 z^6 + \frac{c}{2} z^4 + \frac14(c^2-4) z^2 - z \int x \, \dd\mu(x) - \int x^2 \, \dd\mu(x) - c \, .
\end{displaymath}
We can not find the roots of this polynomial because the two first moments of $\mu$ are unknown. However, we will be able to count its real roots applying Descartes' rule of signs.

\begin{Lemma}[Descartes' rule of signs, see \cite{W}]
Let
\begin{displaymath}
P(X) = a_nX^n + \ldots + a_1X + a_0
\end{displaymath}
be a polynomial with real coefficients. We denote by $p$, resp. $q$, the number of sign changes in the sequence $(a_n,\ldots,a_1,a_0)$, resp. $((-1)^n a_n,\ldots,-a_1,a_0)$, in which we have removed the zeros. Then, the number of positive, resp. negative, roots of $P$ is at most $p$, resp. $q$, and has the same parity as $p$, resp. $q$.
\end{Lemma}

If we distinguish all the possible cases, it easily follows that the polynomial $R$ admits 0, 2, or 4 non-zero real roots, whatever the value of $c \ge -2$ is and whatever the signs of the quantities $\int x \, \dd\mu(x)$ and $\int x^2 \, \dd\mu(x) + c$ are.

In addition to this, every non-real root of $R$ has even multiplicity by Proposition \ref{lemma2} (ii). Since $R$ admits 6 roots, it follows that the multiplicity of 0 is necessarily even.
\begin{itemize}
\item If 0 is not a root of $R$, then $R$ admits at most 4 real roots, thus at least two conjugate non-real roots. But, by Proposition \ref{lemma2} (ii), every non-real root is at least a double root, thus $R$ has in fact at most two real roots. By Proposition \ref{lemma2} (iii), $\mu$ has connected support in this case.
\item If 0 is a root of $R$, then it is at least a double root. Thus $R$ is explicit and we have $R(z) = \frac14 z^2(z^2+c+2)(z^2+c-2)$. This is impossible for $c>-2$ by Proposition \ref{lemma2} (ii). For $c=-2$, this leads to $R(z) = \frac14 z^4(z-2)(z+2)$, thus by Proposition \ref{lemma2} (iii), the support of $\mu$ is $[-2,0]$, $[0,2]$, or $[-2,2]$.
\end{itemize}
In both cases, we have shown that $\mu$ has connected support. \qed
\end{Proof}

\subsection{Singular integral equations}

Euler-Lagrange equations are singular integral equations that can be solved once we know the support of the solution, or at least its number of connected components, thanks to the following result. For a slightly different approach, see \citep{Tri}.

\begin{Theorem}[see {\cite[\S 88]{Musk}}] \label{thmMusk}
Let $L$ be a finite union of intervals $\bigcup_{j=1}^p [a_{2j-1},a_{2j}]$ and let $f$ be a given H\"older continuous function on $L$. The singular integral equation
\begin{displaymath}
\forall x \in L, \quad \fint_L \frac{\varphi(t)}{t-x} \, \dd t = f(x)
\end{displaymath}
admits a H\"older continuous, bounded solution $\varphi$ if and only if $f$ satisfies the following $p$ conditions:
\begin{displaymath}
\forall k \in \llbracket 0,p-1 \rrbracket , \quad \int_L \frac{t^k f(t)}{\prod_{j=1}^{2p} \sqrt{|t-a_j|}} \, \dd t = 0 \, .
\end{displaymath}
In this case, the solution is unique and it is given by
\begin{displaymath}
\forall x \in L, \quad \varphi(x) = -\frac{1}{\pi^2} \prod_{j=1}^{2p} \sqrt{|x-a_j|} \fint_L \frac{f(t)}{(t-x) \prod_{j=1}^{2p} \sqrt{|t-a_j|}} \, \dd t \, .
\end{displaymath}
\end{Theorem}

Applying this theorem to  the quartic potential \eqref{V}, we get the following.

\begin{Proposition} \label{lemma3}
For the potential $V(x) = \frac14 x^4 + \frac{c}{2} x^2$ with $c \ge -2$, the only stationary probability measure with bounded density and connected support is the equilibrium measure $\mu_V$, which is defined by \eqref{mu_V_c>-2} .
\end{Proposition}

\begin{Proof}
Let $\mu$ be a stationary probability measure with bounded density, denoted by $\rho$, and with connected support, denoted by $[a,b]$. By Theorem \ref{thmMusk} applied to $f(x) = -\frac12 V'(x)$ and $p=1$, the existence of $\mu$ is ensured by the condition
\begin{equation} \label{a,b_eq0}
\int_a^b \frac{t^3+ct}{\sqrt{(t-a)(b-t)}} \, \dd t = 0 \, .
\end{equation}
An elementary computation leads to
\begin{displaymath}
\int_a^b \frac{t^3+ct}{\sqrt{(t-a)(b-t)}} \, \dd t = \frac{\pi}{16}(5b^3+3ab^2+3a^2b+5a^3) + c\frac{\pi}{2}(a+b) \, ,
\end{displaymath}
thus condition \eqref{a,b_eq0} reads
\begin{equation} \label{a,b_eq1}
(a+b)(5b^2-2ab+5a^2+8c) = 0 \, .
\end{equation}

Moreover, by Theorem \ref{thmMusk} again, the density of $\mu$ is given by
\begin{eqnarray} \label{mu_eq1}
\rho(x) & = & \frac{\sqrt{(x-a)(b-x)}}{2\pi^2} \fint_a^b \frac{t^3+ct}{(t-x) \sqrt{(t-a)(b-t)}} \, \dd t \nonumber \\
& = & \frac{1}{2\pi} \sqrt{(x-a)(b-x)} \left( x^2 + \frac{a+b}{2} x + \frac38 b^2 + \frac14 ab + \frac38 a^2 + c \right) \, .
\end{eqnarray}
This result has been obtained by standard integral computations. By integrating this expression between $a$ and $b$, since $\rho$ is a probability density function, we get a new constraint on $a$ and $b$:
\begin{equation} \label{a,b_eq2}
\frac{(b-a)^2}{256} (15a^2+18ab+15b^2+16c) = 1 \, .
\end{equation}

The two equations \eqref{a,b_eq1} and \eqref{a,b_eq2} allow us to determine $a$ and $b$. First, Equation \eqref{a,b_eq1} gives three families of possible solutions:
\begin{displaymath}
a = -b , \qquad a = \frac15 b + \frac25 \sqrt{-10c-6b^2} , \qquad a = \frac15 b - \frac25 \sqrt{-10c-6b^2} \, .
\end{displaymath}
Equation \eqref{a,b_eq2} then eliminates some cases. Note first that, if $c$ is non-negative, only the first case would be possible, and that the same situation occurs when $c$ is negative but $b^2 > -\frac53 c$.\\

\noindent $\bullet$ Case 1: $a=-b$.\\
In this case, Equation \eqref{a,b_eq2} gives
\begin{displaymath}
b = \sqrt{\frac23 \left( \sqrt{c^2+12} - c \right)} \, ,
\end{displaymath}
so the density given by \eqref{mu_eq1} becomes
\begin{displaymath}
\rho(x) = \frac{1}{2\pi} \sqrt{b^2-x^2} \left( x^2 + \frac23 c + \frac13 \sqrt{c^2+12} \right) \, .
\end{displaymath}
This is exactly the equilibrium measure of $V$ for $c \ge -2$, see \eqref{mu_V_c>-2}.\\

\noindent $\bullet$ Case 2: $a = \dfrac15 b + \dfrac25 \sqrt{-10c-6b^2}$.\\
Equation \eqref{a,b_eq2} now implies that
\begin{displaymath}
45b^8 + 156cb^6 + (182c^2-552)b^4 + (76c^3-880c)b^2 + 5c^4 - 200c^2 + 2000 = 0 \, .
\end{displaymath}
We will show this is not possible under the conditions $-2 \le c \le 0$ and $0 \le b^2 \le -\frac53 c$. Indeed, we can study the polynomial function
\begin{displaymath}
f : (x,c) \mapsto 45x^4 + 156cx^3 + (182c^2-552)x^2 + (76c^3-880c)x + 5c^4 - 200c^2 + 2000
\end{displaymath}
on the compact set
\begin{displaymath}
K = \left\{ (x,c) \in \R^2 \ | \ -2 \le c \le 0, \ 0 \le x \le -\frac53 c \right\} \, .
\end{displaymath}

The resolution of $\frac{\partial f}{\partial x}(x,c) = \frac{\partial f}{\partial c}(x,c) = 0$ shows that the only critical point of $f$ in $K$ is $(0,0)$. Consequently, $f$ attains its minimum on the boundary of $K$. The study of the three functions
\begin{displaymath}
c \mapsto f(0,c) = 5(c^2-20)^2 \, ,
\end{displaymath}
\begin{displaymath}
x \mapsto f(x,-2) = 45x^4 - 312x^3 + 176x^2 + 1152x + 1280 \, ,
\end{displaymath}
and
\begin{displaymath}
c \mapsto f \left( -\frac53 c, c \right) = \frac{80}{9} (c^2-15)^2
\end{displaymath}
allows us to conclude that the minimum of $f$ on $K$ is attained at $\left( \frac{10}{3}, -2 \right)$ and is equal to $\frac{9680}{9}$. Consequently, $f$ does not vanish on $K$ and Case 2 does not lead to a suitable solution $\mu$.\\

\noindent $\bullet$ Case 3: $a = \dfrac15 b - \dfrac25 \sqrt{-10c-6b^2}$.\\
Very similar computations lead to the fact that the same function $f$ must vanish on the same compact $K$, and thus to the same conclusion.\\

Finally, the only stationary probability measure with bounded density and connected support is indeed the equilibrium measure $\mu_V$. \qed
\end{Proof}

\begin{Remark}
The previous calculations also show that there does not exist a stationary probability measure with bounded density and connected support when $-\sqrt{15} < c < -2$ because, in this situation, Case 1 of the proof leads to a density taking negative values, and Cases 2 and 3 still lead to unsuitable solutions.

In addition to this, the same technique allows us to prove that, when $c<-2$, the only symmetric stationary probability measure having a bounded density and a support with two connected components is the equilibrium measure. We recall that, by Proposition \ref{lemma2_quartic}, there does not exist such a symmetric stationary measure when $c \ge -2$.

Finally, we mention that, for $c=-\sqrt{15}$, there exist two stationary measures with bounded density and connected support. The first one is given by the density
\begin{displaymath}
x \mapsto \frac{1}{2\pi} \sqrt{\left( x - \frac{1}{\sqrt[4]{15}} \right) \left( \frac{5}{\sqrt[4]{15}} - x \right)} \left( x + \frac{4}{\sqrt[4]{15}} \right) \left( x - \frac{1}{\sqrt[4]{15}} \right)
\end{displaymath}
on the interval $\left[ \frac{1}{\sqrt[4]{15}} ; \frac{5}{\sqrt[4]{15}} \right]$ and the second one is its symmetrical measure with respect to the origin. We refer to \cite[Chapter 7]{G} for the detailed computations.
\end{Remark}

\section{Proof of Theorem \ref{theorem}} \label{section_proof}

We are now able to prove Theorem \ref{theorem}. The ideas are as follows. First, thanks to properties of free diffusions stated in Section \ref{section_properties_freeSDE}, we find an accumulation point of $(\mu_t)_{t \ge 0}$ which is a stationary measure with a bounded density. By Propositions \ref{lemma2_quartic} and \ref{lemma3}, this accumulation point is necessarily the equilibrium measure $\mu_V$. We then prove that all accumulation points have the same entropy, using again the estimates of Proposition \ref{lemma1}. Since $\mu_V$ is the unique minimizer of $\Sigma_V$, this proves it is the only accumulation point. A compactness argument allows us to prove the convergence of $(\mu_t)_{t \ge 0}$ towards $\mu_V$.

We emphasize our proof depends on the special potential \eqref{V} only through Propositions \ref{lemma2_quartic} and \ref{lemma3}; the other arguments given below are valid for every potential $V$ satisfying the assumptions of Proposition \ref{lemma1}.\\

From now, for every $t \ge 0$, we denote by $\rho_t$ the density of $\mu_t$.

\begin{Proof}[of Theorem \ref{theorem}]
By \eqref{derivativeSigmaV}, the function $t \mapsto \Sigma_V(\mu_t)$ is decreasing on $[0,+\infty)$. As it is also bounded below (by $\Sigma_V(\mu_V)$), this function admits a finite limit as $t$ goes to infinity. Therefore, by \eqref{derivativeSigmaV} again, there exists a sequence $(t_k)_{k \in \N}$ such that $t_k \to \infty$ and $\frac{\dd}{\dd t} \Sigma_V(\mu_{t_k}) \to 0$ when $k \to \infty$.\\

By Proposition \ref{lemma1} (iii), extracting a further subsequence if necessary, we can assume that the densities $\rho_{t_k}$ converge in the $L^2$-topology to a limit $\rho$. As the $\rho_{t_k}$'s are defined on the compact set $[-M,M]$, $L^2$-convergence implies that $\int \rho_{t_k}(x) \, \dd x$ converges to $\int \rho(x) \, \dd x$, hence the limit $\rho$ is a probability density function defined on $[-M,M]$. We denote by $\mu$ the probability measure associated to $\rho$. By Scheff\'e's lemma, $\mu_{t_k}$ also converges in distribution towards $\mu$.\\

We will now prove that $\mu$ is a stationary probability measure with a bounded density. First, as the densities $\rho_{t_k}$ converge in $L^2([-M,M])$, extracting a further subsequence if necessary, we can assume that they converge almost everywhere on $[-M,M]$. Thus, we have $\| \rho \|_{\infty} \le K_2$.

Furthermore, for all $k \in \N$, we can decompose
\begin{multline} \label{mu_stat_eq1}
\lefteqn{   \left| \int \left| H\mu_{t_k} - \frac12 V' \right|^2 \, \dd\mu_{t_k} - \int \left| H\mu - \frac12 V' \right|^2 \, \dd\mu \right|   } \\
\le \left| \int \left| H\mu_{t_k} - \frac12 V' \right|^2 \, \dd\mu_{t_k} - \int \left| H\mu - \frac12 V' \right|^2 \, \dd\mu_{t_k} \right| \\
+ \left| \int \left| H\mu - \frac12 V' \right|^2 \, \dd\mu_{t_k} - \int \left| H\mu - \frac12 V' \right|^2 \, \dd\mu \right|
\end{multline}
where the integrals are taken over $[-M,M]$.\\
Let us show why the first term in the right-hand side goes to 0 as $k \to +\infty$. Denoting by $K$ a uniform bound on the densities and using the Cauchy-Schwarz inequality, we have for all $k \in \N^*$,
\begin{eqnarray*}
\lefteqn{   \left| \int \left| H\mu_{t_k} - \frac12 V' \right|^2 \, \dd\mu_{t_k} - \int \left| H\mu - \frac12 V' \right|^2 \, \dd\mu_{t_k} \right|   } \\
& \le & K \left| \int \left| H\mu_{t_k}(x) - H\mu(x) \right| \left( |H\mu_{t_k}(x)| + |H\mu(x)| + |V'(x)| \right) \, \dd x \right| \\
& \le & K \| H\mu_{t_k} - H\mu \|_2 \times \left( \| H\mu_{t_k} \|_2 + \| H\mu \|_2 + \| V' \|_2 \right)
\end{eqnarray*}
and by the continuity of the Hilbert transform from $L^2(\R)$ to $L^2(\R)$, $H\mu_{t_k}$ converges to $H\mu$ in $L^2$.\\
On the other hand, it follows from similar arguments and from $\rho \in L^4$ that the second term in the right-hand side of \eqref{mu_stat_eq1} also tends to 0 when $k \to +\infty$. By \eqref{derivativeSigmaV}, we finally have
\begin{displaymath}
0 = \lim_{k \to +\infty} \frac{\dd}{\dd t} \Sigma_V(\mu_{t_k}) = \lim_{k \to +\infty} -2\int \left| H\mu_{t_k} - \frac12 V' \right|^2 \, \dd\mu_{t_k} = -2\int \left| H\mu - \frac12 V' \right|^2 \, \dd\mu \, .
\end{displaymath}
The limit measure $\mu$ is thus a stationary probability measure.\\

As $\mu$ is also a critical measure supported in $\mathbb{R}$, when $V$ is the quartic potential \eqref{V}, by Proposition \ref{lemma2_quartic}, $\mu$ has connected support and, by Proposition \ref{lemma3}, $\mu$ is the equilibrium measure $\mu_V$.\\

We will now prove that two accumulation points for $(\rho_t)_{t \ge 0}$ in the $L^2$-topology must have the same entropy.

Indeed, let $(\rho_{s_k})_{k \in \N}$ be a convergent subsequence from $(\rho_t)_{t \ge 0}$ in the $L^2$-topology. We denote by $\rho$ its limit. Then $\rho$ is a probability density function supported in $[-M,M]$ and it is bounded by $K_1+K_2$, as these properties hold for $\rho_t$ with $t \ge 1$. We denote by $\mu$ the associated probability measure.

By the Cauchy-Schwarz inequality, we have
\begin{displaymath}
\left| \int V(x) \, \dd\mu_{s_k}(x) - \int V(x) \, \dd\mu(x) \right| \le \| \rho_{s_k} - \rho \|_2 \|V\|_2
\end{displaymath}
and
\begin{eqnarray*}
\lefteqn{   \left| \iint \log|x-y| \rho_{s_k}(x) \rho_{s_k}(y) \, \dd x \dd y - \iint \log|x-y| \rho(x) \rho(y) \, \dd x \dd y \right|   } \\
& \le & \left| \iint \log|x-y| \rho_{s_k}(x) (\rho_{s_k}(y) - \rho(y)) \, \dd x \dd y \right| \\
& & \quad + \left| \iint \log|x-y| \rho(y) (\rho_{s_k}(x) - \rho(x)) \, \dd x \dd y \right| \\
& \le & 2(K_1+K_2).\sqrt{2M} \| \rho_{s_k} - \rho \|_2 \left( \iint \log^2|x-y| \, \dd x \dd y \right)^{1/2}
\end{eqnarray*}
for $k$ large enough. Therefore, we get
\begin{displaymath}
\lim_{k \to +\infty} \Sigma_V(\mu_{s_k}) = \Sigma_V(\mu) \, .
\end{displaymath}
This proves the "continuity" of entropy along a solution.

Since the function $t \mapsto \Sigma_V(\mu_t)$ is decreasing, we conclude that two accumulation points lead to the same entropy.\\

This allows us to complete the proof. Indeed, since we proved that the density $\rho_V$ of $\mu_V$ is an accumulation point of $(\rho_t)_{t \ge 0}$ in the $L^2$-topology, since $\mu_V$ is the unique minimizer of free entropy $\Sigma_V$, and since all accumulation points have the same entropy, the only possible accumulation point in the $L^2$-topology is $\rho_V$. But, by Proposition \ref{lemma1} (iii), the $\rho_t$'s, $t \ge 1$, are contained in a compact set $\AAA$ for this topology, so $\rho_t$ converges towards $\rho_V$ in the $L^2$-topology. As we explained at the beginning of this proof, this implies that $\mu_t$ converges in distribution towards $\mu_V$. Since weak convergence and $W_p$-convergence, $p \in [1,+\infty)$, coincide for distributions on a given compact set, the conclusion of Theorem \ref{theorem} follows. \qed
\end{Proof}

\section{Perspectives} \label{section_perspectives}

Many natural questions follow this work.
\begin{itemize}
\item \emph{The case $c<-2$.} Our result uses the fact that, when $c \ge -2$, we have only one critical measure that can be an accumulation point for $(\mu_t)_{t \ge 0}$ (Propositions \ref{lemma2_quartic} and \ref{lemma3}). When $c<-2$, can we describe the critical measures that are candidates to be accumulation points? For instance, there is no critical measure with bounded density and connected support when $-\sqrt{15} < c < -2$. In this case, is the equilibrium measure $\mu_V$ the only suitable critical measure? Is the convergence of the solution of \eqref{FP} towards $\mu_V$ possible in this case?\\
Besides, the value $c = -\sqrt{15}$ appears as the value under which the existence of unilateral critical measures for the quartic potential becomes possible. This threshold also appears in \cite{BT} in a slightly different context. Are the measures described in \cite{BT} the only critical measures?\\
Finally, when $c$ is very negative, can we describe the basins of attraction associated to each possible limit for the solution of the free Fokker-Planck equation?

\item \emph{Other confining potentials.} We only used the special form of the quartic potential in order to get Propositions \ref{lemma2_quartic} and \ref{lemma3}. Do our methods apply in other cases? For instance, can we change the potential $V$, take a higher degree, or consider higher dimensions?

\item \emph{Non-confining potentials.} Several works deal with non-confining potentials. For instance, \cite{AD} studied a cubic potential, and \cite{BIPZ} considered the quartic potential $V(x) = \frac12 x^2 + \frac{g}{4} x^4$ with $g<0$. For these potentials, once the problems of definitions are solved, we can tackle the problem of long-time behaviour. Can we prove a convergence result for the cubic potential or for the quartic potential with $-\frac{1}{12} < g < 0$, as Biane and Speicher conjectured for the latter?
\end{itemize}

\section*{Acknowledgements}

We thank Guilherme Silva for some explanations he gave to us about quadratic differentials when visiting Lille. We also thank one of the anonymous referees for his precise suggestions. M.M. was partially supported by the Labex CEMPI (ANR-11-LABX-0007-01).

{\small
  \setlength{\bibsep}{.5em}
  \bibliographystyle{abbrvnat}
  \bibliography{article2_biblio}

\begin{thebibliography}{35}
\providecommand{\natexlab}[1]{#1}
\providecommand{\url}[1]{\texttt{#1}}
\expandafter\ifx\csname urlstyle\endcsname\relax
  \providecommand{\doi}[1]{doi: #1}\else
  \providecommand{\doi}{doi: \begingroup \urlstyle{rm}\Url}\fi

\bibitem[Allez and Dumaz(2015)]{AD}
R.~Allez and L.~Dumaz.
\newblock Random matrices in non-confining potentials.
\newblock \emph{J. Stat. Phys.}, 160\penalty0 (3):\penalty0 681--714, 2015.
\newblock ISSN 0022-4715.
\newblock \doi{10.1007/s10955-015-1258-1}.
\newblock URL \url{http://dx.doi.org/10.1007/s10955-015-1258-1}.

\bibitem[Anderson et~al.(2010)Anderson, Guionnet, and Zeitouni]{AGZ}
G.~W. Anderson, A.~Guionnet, and O.~Zeitouni.
\newblock \emph{An introduction to random matrices}, volume 118 of
  \emph{Cambridge Studies in Advanced Mathematics}.
\newblock Cambridge University Press, Cambridge, 2010.
\newblock ISBN 978-0-521-19452-5.

\bibitem[Benedetto et~al.(1997)Benedetto, Caglioti, and Pulvirenti]{BCP}
D.~Benedetto, E.~Caglioti, and M.~Pulvirenti.
\newblock A kinetic equation for granular media.
\newblock \emph{RAIRO Mod\'el. Math. Anal. Num\'er.}, 31\penalty0 (5):\penalty0
  615--641, 1997.
\newblock ISSN 0764-583X.

\bibitem[Benedetto et~al.(1998)Benedetto, Caglioti, Carrillo, and
  Pulvirenti]{BCCP}
D.~Benedetto, E.~Caglioti, J.~A. Carrillo, and M.~Pulvirenti.
\newblock A non-{M}axwellian steady distribution for one-dimensional granular
  media.
\newblock \emph{J. Statist. Phys.}, 91\penalty0 (5-6):\penalty0 979--990, 1998.
\newblock ISSN 0022-4715.
\newblock \doi{10.1023/A:1023032000560}.
\newblock URL \url{http://dx.doi.org/10.1023/A:1023032000560}.

\bibitem[Bertola and Tovbis(2015)]{BT}
M.~Bertola and A.~Tovbis.
\newblock Asymptotics of orthogonal polynomials with complex varying quartic
  weight: global structure, critical point behavior and the first {P}ainlev\'e
  equation.
\newblock \emph{Constr. Approx.}, 41\penalty0 (3):\penalty0 529--587, 2015.
\newblock ISSN 0176-4276.
\newblock \doi{10.1007/s00365-015-9288-0}.
\newblock URL \url{http://dx.doi.org/10.1007/s00365-015-9288-0}.

\bibitem[Biane(1997)]{Biane}
P.~Biane.
\newblock Free {B}rownian motion, free stochastic calculus and random matrices.
\newblock In \emph{Free probability theory ({W}aterloo, {ON}, 1995)}, volume~12
  of \emph{Fields Inst. Commun.}, pages 1--19. Amer. Math. Soc., Providence,
  RI, 1997.

\bibitem[Biane and Speicher(1998)]{BSp1}
P.~Biane and R.~Speicher.
\newblock Stochastic calculus with respect to free {B}rownian motion and
  analysis on {W}igner space.
\newblock \emph{Probab. Theory Related Fields}, 112\penalty0 (3):\penalty0
  373--409, 1998.
\newblock ISSN 0178-8051.
\newblock \doi{10.1007/s004400050194}.
\newblock URL \url{http://dx.doi.org/10.1007/s004400050194}.

\bibitem[Biane and Speicher(2001)]{BS}
P.~Biane and R.~Speicher.
\newblock Free diffusions, free entropy and free {F}isher information.
\newblock \emph{Ann. Inst. H. Poincar\'e Probab. Statist.}, 37\penalty0
  (5):\penalty0 581--606, 2001.
\newblock ISSN 0246-0203.
\newblock \doi{10.1016/S0246-0203(00)01074-8}.
\newblock URL \url{http://dx.doi.org/10.1016/S0246-0203(00)01074-8}.

\bibitem[Bolley et~al.(2010)Bolley, Guillin, and Malrieu]{BGM}
F.~Bolley, A.~Guillin, and F.~Malrieu.
\newblock Trend to equilibrium and particle approximation for a weakly
  selfconsistent {V}lasov-{F}okker-{P}lanck equation.
\newblock \emph{M2AN Math. Model. Numer. Anal.}, 44\penalty0 (5):\penalty0
  867--884, 2010.
\newblock ISSN 0764-583X.
\newblock \doi{10.1051/m2an/2010045}.
\newblock URL \url{http://dx.doi.org/10.1051/m2an/2010045}.

\bibitem[Bolley et~al.(2012)Bolley, Gentil, and Guillin]{BGG1}
F.~Bolley, I.~Gentil, and A.~Guillin.
\newblock Convergence to equilibrium in {W}asserstein distance for
  {F}okker-{P}lanck equations.
\newblock \emph{J. Funct. Anal.}, 263\penalty0 (8):\penalty0 2430--2457, 2012.
\newblock ISSN 0022-1236.
\newblock \doi{10.1016/j.jfa.2012.07.007}.
\newblock URL \url{http://dx.doi.org/10.1016/j.jfa.2012.07.007}.

\bibitem[Bolley et~al.(2013)Bolley, Gentil, and Guillin]{BGG2}
F.~Bolley, I.~Gentil, and A.~Guillin.
\newblock Uniform convergence to equilibrium for granular media.
\newblock \emph{Arch. Ration. Mech. Anal.}, 208\penalty0 (2):\penalty0
  429--445, 2013.
\newblock ISSN 0003-9527.
\newblock \doi{10.1007/s00205-012-0599-z}.
\newblock URL \url{http://dx.doi.org/10.1007/s00205-012-0599-z}.

\bibitem[Br{\'e}zin et~al.(1978)Br{\'e}zin, Itzykson, Parisi, and Zuber]{BIPZ}
E.~Br{\'e}zin, C.~Itzykson, G.~Parisi, and J.~B. Zuber.
\newblock Planar diagrams.
\newblock \emph{Comm. Math. Phys.}, 59\penalty0 (1):\penalty0 35--51, 1978.
\newblock ISSN 0010-3616.

\bibitem[Carrillo et~al.(2003)Carrillo, McCann, and Villani]{CMCV}
J.~A. Carrillo, R.~J. McCann, and C.~Villani.
\newblock Kinetic equilibration rates for granular media and related equations:
  entropy dissipation and mass transportation estimates.
\newblock \emph{Rev. Mat. Iberoamericana}, 19\penalty0 (3):\penalty0 971--1018,
  2003.
\newblock ISSN 0213-2230.
\newblock \doi{10.4171/RMI/376}.
\newblock URL \url{http://dx.doi.org/10.4171/RMI/376}.

\bibitem[Carrillo et~al.(2015)Carrillo, Castorina, and Volzone]{CCV}
J.~A. Carrillo, D.~Castorina, and B.~Volzone.
\newblock Ground states for diffusion dominated free energies with logarithmic
  interaction.
\newblock \emph{SIAM J. Math. Anal.}, 47\penalty0 (1):\penalty0 1--25, 2015.
\newblock ISSN 0036-1410.
\newblock \doi{10.1137/140951588}.
\newblock URL \url{http://dx.doi.org/10.1137/140951588}.

\bibitem[Cattiaux et~al.(2008)Cattiaux, Guillin, and Malrieu]{CGM}
P.~Cattiaux, A.~Guillin, and F.~Malrieu.
\newblock Probabilistic approach for granular media equations in the
  non-uniformly convex case.
\newblock \emph{Probab. Theory Related Fields}, 140\penalty0 (1-2):\penalty0
  19--40, 2008.
\newblock ISSN 0178-8051.
\newblock \doi{10.1007/s00440-007-0056-3}.
\newblock URL \url{http://dx.doi.org/10.1007/s00440-007-0056-3}.

\bibitem[C{\'e}pa and L{\'e}pingle(1997)]{CL}
E.~C{\'e}pa and D.~L{\'e}pingle.
\newblock Diffusing particles with electrostatic repulsion.
\newblock \emph{Probab. Theory Related Fields}, 107\penalty0 (4):\penalty0
  429--449, 1997.
\newblock ISSN 0178-8051.
\newblock \doi{10.1007/s004400050092}.
\newblock URL \url{http://dx.doi.org/10.1007/s004400050092}.

\bibitem[Chan(1992)]{Chan}
T.~Chan.
\newblock The {W}igner semi-circle law and eigenvalues of matrix-valued
  diffusions.
\newblock \emph{Probab. Theory Related Fields}, 93\penalty0 (2):\penalty0
  249--272, 1992.
\newblock ISSN 0178-8051.
\newblock \doi{10.1007/BF01195231}.
\newblock URL \url{http://dx.doi.org/10.1007/BF01195231}.

\bibitem[Demengel and Demengel(2012)]{Demengel}
F.~Demengel and G.~Demengel.
\newblock \emph{Functional spaces for the theory of elliptic partial
  differential equations}.
\newblock Universitext. Springer, London; EDP Sciences, Les Ulis, 2012.
\newblock ISBN 978-1-4471-2806-9; 978-2-7598-0698-0.
\newblock \doi{10.1007/978-1-4471-2807-6}.
\newblock URL \url{http://dx.doi.org/10.1007/978-1-4471-2807-6}.
\newblock Translated from the 2007 French original by Reinie Ern{\'e}.

\bibitem[Dyson(1962)]{Dyson}
F.~J. Dyson.
\newblock A {B}rownian-motion model for the eigenvalues of a random matrix.
\newblock \emph{J. Mathematical Phys.}, 3:\penalty0 1191--1198, 1962.
\newblock ISSN 0022-2488.

\bibitem[Fontbona(2004)]{Font}
J.~Fontbona.
\newblock Uniqueness for a weak nonlinear evolution equation and large
  deviations for diffusing particles with electrostatic repulsion.
\newblock \emph{Stochastic Process. Appl.}, 112\penalty0 (1):\penalty0
  119--144, 2004.
\newblock ISSN 0304-4149.
\newblock \doi{10.1016/j.spa.2004.01.008}.
\newblock URL \url{http://dx.doi.org/10.1016/j.spa.2004.01.008}.

\bibitem[Groux(2016)]{G}
B.~Groux.
\newblock \emph{Grandes d\'eviations de matrices al\'eatoires et \'Equation de
  Fokker-Planck libre}.
\newblock {PhD} thesis, {Universit{\'e} Paris-Saclay}, 2016.
\newblock URL \url{https://tel.archives-ouvertes.fr/tel-01507380}.

\bibitem[Huybrechs et~al.(2014)Huybrechs, Kuijlaars, and Lejon]{HKL}
D.~Huybrechs, A.~B.~J. Kuijlaars, and N.~Lejon.
\newblock Zero distribution of complex orthogonal polynomials with respect to
  exponential weights.
\newblock \emph{J. Approx. Theory}, 184:\penalty0 28--54, 2014.
\newblock ISSN 0021-9045.
\newblock \doi{10.1016/j.jat.2014.05.002}.
\newblock URL \url{http://dx.doi.org/10.1016/j.jat.2014.05.002}.

\bibitem[Johansson(1998)]{J}
K.~Johansson.
\newblock On fluctuations of eigenvalues of random {H}ermitian matrices.
\newblock \emph{Duke Math. J.}, 91\penalty0 (1):\penalty0 151--204, 1998.
\newblock ISSN 0012-7094.
\newblock \doi{10.1215/S0012-7094-98-09108-6}.
\newblock URL \url{http://dx.doi.org/10.1215/S0012-7094-98-09108-6}.

\bibitem[Kuijlaars and Silva(2015)]{KS}
A.~B.~J. Kuijlaars and G.~L.~F. Silva.
\newblock S-curves in polynomial external fields.
\newblock \emph{J. Approx. Theory}, 191:\penalty0 1--37, 2015.
\newblock ISSN 0021-9045.
\newblock \doi{10.1016/j.jat.2014.04.002}.
\newblock URL \url{http://dx.doi.org/10.1016/j.jat.2014.04.002}.

\bibitem[Li et~al.(2014)Li, Li, and Xie]{LLX}
S.~Li, X.~Li, and Y.~Xie.
\newblock On the {L}aw of {L}arge {N}umbers for the empirical measure process
  of generalized {D}yson {B}rownian {M}otion.
\newblock \emph{arXiv:1407.7234v2}, 2014.

\bibitem[Malrieu(2003)]{Malrieu}
F.~Malrieu.
\newblock Convergence to equilibrium for granular media equations and their
  {E}uler schemes.
\newblock \emph{Ann. Appl. Probab.}, 13\penalty0 (2):\penalty0 540--560, 2003.
\newblock ISSN 1050-5164.
\newblock \doi{10.1214/aoap/1050689593}.
\newblock URL \url{http://dx.doi.org/10.1214/aoap/1050689593}.

\bibitem[Mart{\'{\i}}nez-Finkelshtein and Rakhmanov(2011)]{MFR}
A.~Mart{\'{\i}}nez-Finkelshtein and E.~A. Rakhmanov.
\newblock Critical measures, quadratic differentials, and weak limits of zeros
  of {S}tieltjes polynomials.
\newblock \emph{Comm. Math. Phys.}, 302\penalty0 (1):\penalty0 53--111, 2011.
\newblock ISSN 0010-3616.
\newblock \doi{10.1007/s00220-010-1177-6}.
\newblock URL \url{http://dx.doi.org/10.1007/s00220-010-1177-6}.

\bibitem[Muskhelishvili(1972)]{Musk}
N.~I. Muskhelishvili.
\newblock \emph{Singular integral equations}.
\newblock Wolters-Noordhoff Publishing, Groningen, 1972.
\newblock Boundary problems of functions theory and their applications to
  mathematical physics, Revised translation from the Russian, edited by J. R.
  M. Radok, Reprinted.

\bibitem[Rogers and Shi(1993)]{RS}
L.~C.~G. Rogers and Z.~Shi.
\newblock Interacting {B}rownian particles and the {W}igner law.
\newblock \emph{Probab. Theory Related Fields}, 95\penalty0 (4):\penalty0
  555--570, 1993.
\newblock ISSN 0178-8051.
\newblock \doi{10.1007/BF01196734}.
\newblock URL \url{http://dx.doi.org/10.1007/BF01196734}.

\bibitem[Saff and Totik(1997)]{ST}
E.~B. Saff and V.~Totik.
\newblock \emph{Logarithmic potentials with external fields}, volume 316 of
  \emph{Grundlehren der Mathematischen Wissenschaften [Fundamental Principles
  of Mathematical Sciences]}.
\newblock Springer-Verlag, Berlin, 1997.
\newblock ISBN 3-540-57078-0.
\newblock \doi{10.1007/978-3-662-03329-6}.
\newblock URL \url{http://dx.doi.org/10.1007/978-3-662-03329-6}.
\newblock Appendix B by Thomas Bloom.

\bibitem[Tricomi(1957)]{Tri}
F.~G. Tricomi.
\newblock \emph{Integral equations}.
\newblock Pure and Applied Mathematics. Vol. V. Interscience Publishers, Inc.,
  New York; Interscience Publishers Ltd., London, 1957.

\bibitem[Tugaut(2013{\natexlab{a}})]{Tugaut1}
J.~Tugaut.
\newblock Self-stabilizing processes in multi-wells landscape in
  {$\mathbb{R}^d$}-convergence.
\newblock \emph{Stochastic Process. Appl.}, 123\penalty0 (5):\penalty0
  1780--1801, 2013{\natexlab{a}}.
\newblock ISSN 0304-4149.
\newblock \doi{10.1016/j.spa.2012.12.003}.
\newblock URL \url{http://dx.doi.org/10.1016/j.spa.2012.12.003}.

\bibitem[Tugaut(2013{\natexlab{b}})]{Tugaut2}
J.~Tugaut.
\newblock Convergence to the equilibria for self-stabilizing processes in
  double-well landscape.
\newblock \emph{Ann. Probab.}, 41\penalty0 (3A):\penalty0 1427--1460,
  2013{\natexlab{b}}.
\newblock ISSN 0091-1798.
\newblock \doi{10.1214/12-AOP749}.
\newblock URL \url{http://dx.doi.org/10.1214/12-AOP749}.

\bibitem[Villani(2003)]{Vill}
C.~Villani.
\newblock \emph{Topics in optimal transportation}, volume~58 of \emph{Graduate
  Studies in Mathematics}.
\newblock American Mathematical Society, Providence, RI, 2003.
\newblock ISBN 0-8218-3312-X.

\bibitem[Weisstein()]{W}
E.~W. Weisstein.
\newblock Descartes' sign rule.
\newblock \emph{{M}ath{W}orld--{A} {W}olfram {W}eb resource}.
\newblock URL \url{http://mathworld.wolfram.com/DescartesSignRule.html}.

\end{thebibliography}
}

\end{document}